\ifx\selectfont\undefined
\documentstyle[12pt]{article}
\else
\documentstyle[12pt,oldlfont]{article}
\fi

\newcount\hh
\newcount\mm
\mm=\time
\hh=\time
\divide\hh by 60
\divide\mm by 60
\multiply\mm by 60
\mm=-\mm
\advance\mm by \time
\def\hhmm{\number\hh:\ifnum\mm<10{}0\fi\number\mm}
\makeatletter
\@addtoreset{equation}{section}
\makeatother
\newtheorem{fact}{Fact}[section]
\newtheorem{thm}[fact]{Theorem}
\newtheorem{prop}[fact]{Proposition}
\newtheorem{lemma}[fact]{Lemma}
\newtheorem{cor}[fact]{Corollary}

\newbox\nrmbox
\setbox\nrmbox=\hbox{$\Vert$}
\def\nrmrule{\vrule height\ht\nrmbox depth1.1\dp\nrmbox}
\setbox\nrmbox=%
  \hbox{\kern0.14em\nrmrule\kern0.14em\nrmrule\kern0.14em\nrmrule\kern0.14em}
\newcommand{\Snorm}[1]%
  {\copy\nrmbox#1\copy\nrmbox\kern-0.03em\lower.4ex\hbox{}}
\newbox\Bnrmbox
\setbox\Bnrmbox=\hbox{$\Bigl\Vert\Bigr.$}
\def\Bnrmrule{\vrule height1.03\ht\Bnrmbox depth\dp\Bnrmbox}
\setbox\Bnrmbox=%
 \hbox{\kern0.14em\Bnrmrule\kern0.14em\Bnrmrule\kern0.14em\Bnrmrule\kern0.14em}
\newcommand{\BSnorm}[1]%
  {\copy\Bnrmbox#1\copy\Bnrmbox\kern-0.03em\lower.4ex\hbox{}}
\newcommand{\TTnorm}[1]{\|#1\|}
\newcommand{\BTTnorm}[1]{\Bigl\|#1\Bigr\|}
\newcommand{\rem}{\noindent{\bf Remark{\ \ }}}
\newcommand{\proof}{{\noindent\bf Proof{\ \ }}}
\newcommand{\qed}{\mbox{}\hfill\(\Box\)\bigskip}

\newcommand{\Rn}[1]{\mbox{{\it I\kern -0.25emR}$^{\,{#1}}$}}
\newcommand{\Kn}[1]{\mbox{{\it I\kern -0.25emK}$^{\,{#1}}$}}

\newcommand{\NN}{\mbox{{\it I\kern -0.25emN}}}
\newcommand{\al}{\alpha}
\newcommand{\bt}{\beta}
\newcommand{\de}{\delta}
\newcommand{\De}{\Delta}
\newcommand{\ep}{\varepsilon}

\newcommand{\la}{\lambda}
\newcommand{\La}{\Lambda}
\newcommand{\ga}{\gamma}
\newcommand{\spn}{\mathop{\rm span\,}}

\newcommand{\lust}{\mathop{\rm \, l.u.st\,}}

\newcommand{\ie}{{\em i.e.,\/}}
\newcommand{\eg}{{\em e.g.,\/}}
\newcommand{\cf}{{\em cf.\/}}
\newcommand{\restr}[2]{#1 \,\vert\sb {#2}}
\title{Banach spaces without local unconditional structure\\
}
\author{Ryszard A. Komorowski
\and
Nicole Tomczak-Jaegermann}
\date{\ }
\begin{document}
\maketitle

\begin{abstract}
  For a large class of Banach spaces, a general construction
  of subspaces without local unconditional structure is presented.  As
  an application it is shown that every Banach space of finite cotype
  contains either $l_2$ or a subspace without unconditional basis,
  which admits a Schauder basis.
  Some other interesting applications and corollaries follow.
\end{abstract}
%


\section{Introduction}

In this paper we present, for a large class of
Banach spaces,  a general construction of
subspaces with a basis which have no
local unconditional structure.
The method works for a direct sum of several
Banach spaces with bases which have certain unconditional
properties. It is then  applied
to Banach spaces with unconditional
basis, to show  that if such a space $X$ is of finite cotype
and it does not contain an isomorphic copy of $l_2$,
then $X$ contains a subspace with a basis and without
local unconditional  structure.
As an immediate consequence we get that
if all subspaces of
a Banach space $X$ have unconditional basis
then $X$ is $l_2$ saturated (\ie\  every infinite-dimensional
subspace of $X$ contains a copy of $l_2$).
In particular, if $X$ is a homogeneous Banach space
non-isomorphic to a Hilbert space
(\ie\   $X$ is isomorphic to its every infinite-dimensional
subspace) then
$X$ must not have an unconditional basic
sequence.

We also discuss several other situations.
Let us only  mention here that our method provides a
uniform construction of  subspaces without local
unconditional structure which still
have Gordon--Lewis property in all $L_p$ spaces for $1 \le p < \infty$,
$p \ne 2$, and in all $p$-convexified Tsirelson spaces and their duals
$1 \le p < \infty$.

The technique developed here is based on the approach first
introduced by W.~B.~Johnson, J.~Lindenstrauss and G.~Schechtman
in [J-L-S] for investigating the Kalton--Peck space, which was the
first example of a Banach space which admits 2-dimensional
unconditional decomposition but has no unconditional basis.
This approach was  refined by T.~Ketonen in [Ke]
and subsequently generalized by A.~Borzyszkowski in [B],
for subspaces of $L_p$, with $ 1 \le p <2$.

The essential idea of the approach from  [J-L-S], [Ke] and [B]
is summarized (and slightly generalized for our purpose)
in  Section 1. In the same section we also introduce all
definitions and notations.
Our general construction is presented in Section 2. The
additional ingredient which appears here consists
of  an ordered sequence of partitions of 
natural numbers, which allows to replace some ``global''
arguments used before by ``local'' analogues.
In Section 3  we prove the main application on subspaces
of spaces with an unconditional basis.  Other applications
and corollaries are discussed in Section~4.

\smallskip
After this paper was sent for publication we learnt
about a spectacular structural theorem
just proved by  W.~T.~Gowers.
This theorem combined with our
Theorem 4.2 and a result from [G-M] shows
that a homogeneous Banach space is isomorphic
to a Hilbert space, thus solving in
the positive the so-called homogeneous
space problem. More details can be found
in the  paper by Gowers [G].

\smallskip
The contribution of the first named author
is a part of his Ph.~D.~Thesis written at the University of
Alberta under a supervision of  the second named author.
During the final work on the paper the first named author
was supported by KBN.

\smallskip
\section{Notation and preliminaries}

We use the standard notation from the Banach space theory,
which can be found \eg\ in [L-T.1], [L-T.2] and [T],
together with  all  terminology  not explained here.
In particular, the fundamantal concepts of a basis
and a Schauder decomposition can be found in [L-T.1],
1.a.1 and 1.g.1, respectively.

Let us only recall fundamental notions related
to unconditionality.

A basis $\{e_j\}_j$ in a Banach space $X$ is called
unconditional, if there is a constant $C$
such that for every $x = \sum_j t_j e_j \in X$
one has
$\|\sum_j \ep_j t_j e_j \| \le C \|x\|$,
for all $\ep_j = \pm 1$ for $j=1, 2, \ldots$.
The infimum of constants $C$ is denoted by
$\mbox{unc\,}(\{e_j\})$.
The basis is called 1-unconditional, if
$\mbox{unc\,}(\{e_j\}) = 1$.

A Schauder decomposition  $\{Z_{k}\}_k$
of a Banach space  $X$
is called $C$-un\-con\-di\-tio\-nal,
for some constant $C$, if for all finite
sequences $\{z_k\}$ with $z_k \in Z_k$ for all $k$,
one has
$\|\sum_k \ep_k z_k \| \le C \|\sum_k z_k \|$.
For a subset $K \subset \NN$,
by $Y_K$ denote
$\overline{\spn}\,   [Z_k] _{k \in K}$.

A Banach space $X$ has  local unconditional
structure if there is $C \ge 1$ such
for every  finite-dimensional subspace $X_0 \subset X$
there exist a Banach space $F $ with
a 1-unconditional basis  and operators
$u_0: X_0 \to F$ and $w_0: F \to X$
such that the natural embedding
$j:  X_0 \to X$ admits  a factorization $j = w_0\, u_0$
and $\|u_0\|\, \|w_0\| \le C$.  The infimum of constants $C$
is denoted by $ \lust  (X)$.

\medskip

We will also use several more specific notation.
Let  $F$  be a Banach space with a
basis  $\{f_{l}\}_{l}$. For a subset $A \subset \NN$,
by $\restr{F}{A}$ we denote $\overline{\spn} [f_l]_{l \in A}$.
If  $F'$  is another space with a
basis  $\{f_{l}'\}_{l}$,
by $I: F \to F'$ we denote the formal
identity operator, \ie\
$I(x)= \sum_l t_l f_l'$, for $x = \sum_l t_l f_l \in F$.
With some abuse of  notation, we will occasionally write
$\|I: F \to F'\|=\infty$ when this operator
is not bounded.

We say that a basis $\{f_{l}\}_{l}$ dominates
(resp. is dominated by) $\{f_{l}'\}_{l}$,
if the operator $I: F \to F'$
(resp.   $I: F' \to F$) is bounded.
If the bases in $F$ and $F'$ are fixed and they
are equivalent, by $de (F, F')$ we denote
the equivalence constant,
\begin{equation}
de (F, F') = \|I: F \to F'\|\,  \|I : F' \to F\|;
  \label{def_de}
\end{equation}
and we set  $de (F, F')= \infty$  if the bases
are not equivalent.

By  $D(F\oplus F')$ we denote
the diagonal subspace of $F\oplus F'$,
\ie\ the subspace with the basis
$\{(f_{ j}+ f_{j}')/ \|f_{ j} + f_{j}'\|\}_j$;
an analogous notation will be also used
for a larger (but finite) number of terms.

\medskip
The following  proposition is a version
of a fundamental criterium due to Ketonen [Ke]
and Borzyszkowski [B]. Since a modification
of original arguments would be rather messy,
we provide a shorter direct  proof.

\begin{prop}
\label{borz}
Let $Y$ be a Banach space of cotype $r$,
for some $r < \infty$,
which has a Schauder decomposition $\{Z_{k}\}_k$,
with $\dim Z_k =2$, for $k=1, 2, \ldots$.
If $Y$ has local  unconditional structure  then
there exists a linear, not neccessarily bounded,
operator
$T:   \spn [ Z_k]_k  \to  \spn [ Z_k] _k $
such that
\begin{description}
\item[(i)] $T(Z_k)\subset Z_k$ for $k=1, 2, \ldots$;
\item[(ii)]
If, for some $K \subset \NN$ and
 some  $C \ge 1$,
the decomposition
$\{Z_{k}\}_{ k \in K}$ of $Y_K$
is  $C$-unconditional, then
\begin{equation}
\|\restr{T}{Y_K}: Y_K \to Y_K \| \le  C^2
       \psi  \lust (Y),
  \label{borz_ii}
\end{equation}
where $\psi = \psi( r, C_r(Y))$
depends on  $r$ and the cotype $r$ constant
$ C_r(Y) $  of $Y$ only;
\item[(iii)] $\inf_{\la }\|\restr{T}{Z_k}- \la I_{Z_k}\|
    \ge 1/8$,  for   $k=1, 2, \ldots$.
\end{description}
\end{prop}

The proof  requires a  fact already used in a more
general form in [B]. For sake of completeness
and clarity  of the exposition, we sketch
the proof here.
\begin{lemma}
  \label{borz_cotype}
Let $Y$ be a Banach space of cotype $r$  which has
local  unconditional structure, and let $q > r$.
For every $\ep >0$ and
every  finite-dimensional subspace $Y_0 \subset Y$
there exist a Banach space $E $ with
a 1-unconditional basis  which is
$q$-concave, and operators
$u: Y_0 \to E$ and $w: E \to Y$
such that the natural embedding
$j:  Y_0 \to Y$ admits  a factorization $j = w\, u$
and $\|u\|\, \|w\| \le (1+ \ep) \lust  (Y)$.
Moreover, the $q$-concavity constant of $E$
satisfies $M_{(q)}(E)\le \phi $
where $\phi = \phi  (r, q, C_r(Y))$
depends on $r$, $q$ and the cotype $r$ constant of $Y$
only.
\end{lemma}
\proof
Given $\ep >0$ and $Y_0 $, let
$F$ be a space with a 1-unconditional basis
$\{f_i\}_i$
and let $u_0: Y_0 \to F$ and $w_0: F \to Y$
be  such that $j = w_0 u_0$ and
$\|w_0\|\, \|u_0\| \le (1+\ep) \lust (Y)$.
It can be clearly assumed that $F$ is finite-dimensional,
say $\dim F = N$.
Let  $\{f_i^*\}_i$ be the biorthogonal functionals.

We let $E$ to be $\Rn{N}$ with the norm $\|\cdot\|_E$
defined by
$$
\|(t_i)_i\|_E = \sup_{\ep_i = \pm 1}
   \|w_0 \Bigl(\sum_i \ep_i t_i f_i \Bigr)\|
\qquad \mbox{for\ } (t_i) \in \Rn{N}.
$$
We also set,
$ u(x) = \Bigl( f_i^*(u_0 x)\Bigr)_i $,
for $x \in Y_0$ and
$w \Bigl( (t_i)_i\Bigr) = \sum_i t_i w_0 f_i$,
for $(t_i) \in E$.

It is easy to check that $w u (x) = x$,
for $x \in Y_0$ and that $\|u\| \le \|w_0\|\, \|u_0\|$
and $\|w\| = 1$.
Clearly, the standard unit vector basis is
1-unconditional in $E$. Using the cotype $r$
of $Y$, it can be checked that $E$ satisfies
a lower $r$ estimate with the constant $C_r(Y)$.
Thus $E$ is  $q$-concave for every $q > r$
with the $q$-concavity constant $M_{(q)}(E)$
depending on $q$, $r$ and $C_r(Y)$.
(\cf\ [L-T.2] 1.f.7).
\qed

\noindent{\bf Proof of Proposition~\ref{borz}\ }
Assume that $Y$ has the local  unconditional structure.
It is enough to construct a sequence  of operators
$T_n :  \spn [ Z_k]_k  \to  \spn [ Z_k] _k $, such
that for every $n$, the operator $T_n$
satisfies (i), (ii) and
\begin{description}
\item[(iii')] $\inf_{\la }\|\restr{T_n}{Z_k}- \la I_{Z_k}\|\ge 1/8$,
for   $k=1, 2, \ldots, n$.
\end{description}
Then the existence of the operator $T$ will follow
by Cantor's diagonal procedure and Banach--Steinhaus theorem.

Fix $n$ and $\ep >0$, set $q = 2r$. Let $E$ with
a 1-unconditional basis $\{e_j\}_j$ and operators
$u: Y_{\{1, \ldots, n\}} \to E$ and $w: E \to Y$
be given by Lemma~\ref{borz_cotype},
such that  $j = w\, u$
and $\|u\|\, \|w\| \le (1+\ep) \lust (Y)$;
moreover, $E$ is $2r$-concave.

Let $P_k: Y \to Z_k$ be the natural projection
onto $Z_k$, for $k=1, 2, \ldots$.
For a sequence of signs $\Theta =\{\theta_j\}$,
with $\theta_j = \pm 1$ for $j=1, 2, \ldots$,
define $\La_{\Theta}: E \to E$ by
$\La_{\Theta}(y) = \sum_j \theta_j t_j e_j $,
for $y = \sum_j t_j e_j \in E $.
Then  $\|\La_{\Theta}\|=1$.

For every $k=1, 2, \ldots$ pick a sequence of signs
$\Theta_k$ such that
$$
\sup_\Theta \inf_\la \|P_k w\La_\Theta u P_k - \la {I}_{Z_k}\|
\le (4/3) \inf_\la \|P_k w\La_{\Theta_k}u P_k - \la {I}_{Z_k}\|.
$$

Define $T_n: \spn [ Z_k]_k  \to  \spn [ Z_k]_k $ by
$$
  T_n(y) = \sum_{k=1}^n P_k w  \La_{\Theta_k} u P_k (y)
\qquad \mbox{for\ } y = \sum_k z_k \in \spn [ Z_k]_k.
$$

Clearly (i) follows just from the definition  of $T_n$.
To prove (ii),
let $K_n = K \cap \{1, \ldots, n\}$.
Let $r_k$ denote the Rademacher functions on $[0,1]$.
Since  $(E,\TTnorm{\cdot})$ is a $2r$-concave
Banach lattice with the  $2r$-concavity
constant depending on $r$ and $C_r(Y)$,
and also the decomposition $\{Z_{k}\}_{ k \in K}$ of $Y_K$ is
$C$-unconditional, by Khintchine--Maurey's inequality
(\cf\ \eg\ [L-T.2], 1.d.6) we have,
for $y \in Y_{K_n}$,
\begin{eqnarray*}
  \|\restr{T_n}{Y_K}(y)\| & = &
  \Bigl\|\sum_{k=1}^n P_k w  \La_{\Theta_k}u P_k (y)\Bigr\|
  = \Bigl\|\sum_{k\in {K_n}} P_k w  \La_{\Theta_k}u P_k (y)\Bigr\| \\
  & = & \Bigl\|\int_0^1\Bigl(\sum_{k\in {K_n}} r_k(t) P_k \Bigr)
  \Bigl(\sum_{k\in {K_n}} r_k (t) w\La_{\Theta_k}u P_k (y) \Bigr) dt\Bigr\| \\
  & \le & \sup_{0 \le t \le 1}
  \Bigl\|\sum_{k\in {K_n}} r_k (t) P_k \Bigr\|
      \|w\|
  \int_0^1 \BTTnorm{ \sum_{k\in {K_n}} r_k (t) \La_{\Theta_k} u P_k (y) } dt \\
  & \le & C\, M \|w\| \BTTnorm{
    \Bigl( \sum_{k\in {K_n}}|\La_{\Theta_k}u P_k (y)|^2\Bigr)^{1/2}}\\
  & = & C\, M \|w\| \BTTnorm{
    \Bigl( \sum_{k\in {K_n}}| u P_k (y)|^2\Bigr)^{1/2}}\\
  & \le &  C\, M^2  \|w\| \int_0^1 \BTTnorm{
     \sum_{k\in {K_n}} r_k (t) u  P_k (y)}  dt \\
  & \le &   C\, M^2 \|w\|\,\|u\| \int_0^1
  \Bigl\|  \sum_{k\in {K_n}} r_k (t)   P_k (y) \Bigr\| dt \\
  & \le &  C^2 M^2\, (1 + \ep) \lust (Y)\, \|y\|.
\end{eqnarray*}
The constant $M$ depends on $r$ and
$M_{(2r)}(E)$, hence, implicitely,
on $r$ and $C_r(Y)$; so the function $\psi$
so obtained satisfies the requirements
of (ii).

To prove  (iii'),
fix an arbitrary $k=1, 2, \ldots, n$.
Consider the 4-dimensional space $H$ of all linear
operators on $Z_k$ and the subspace $H_0= \spn [I_{Z_k}]$
spanned by the identity operator on $Z_k$.
Consider the quotient space $H / H_0$ and
for $R \in H$, let $\widetilde{R}$ be the
canonical image of $R$ in  $H / H_0$.

Denote the biorthogonal functionals
to the basis $\{e_j\}_j$ in $E$ by  $\{e_j^*\}_j$
and consider   operators
$R_j= P_k w (e_j^* \otimes e_j)u P_k$
on $Z_k$. Since $\dim R_j(Z_k)=1 < 2$,
it is easy to see that
for every $j=1, 2,\ldots$, one has
$$
\|\widetilde{R}_j\|= \inf_\la \|R_j - \la I_{Z_k}\| \ge (1/2)\|R_j\|.
$$
%

Also recall that if $F$ is an  $m$-dimensional space
then for any vectors $\{x_j\}_j$ in $F$ one has
$$
\sup_{\theta_j= \pm 1} \Bigl\|\sum_j \theta_j x_j\Bigr\|
\ge (1/m) \sum_j \| x_j\|.
$$
This is a restatement of the estimate for the 1-summing
norm of the identity on $F$, $\pi_1 (I_F)\le m$, and it is
a simple consequence of the Auerbach lemma (\cf\ \eg\ [T]).

So by the definition of $T_n$ and by
the choice of $\Theta_k$ and the above estimates we get
\begin{eqnarray*}
\lefteqn{  \inf_{\la }\|\restr{T_n}{Z_k}- \la I_{Z_k}\|
\ge  (3/4) \sup_\Theta \inf_\la \|P_k w\La_\Theta uP_k - \la {I}_{Z_k}\|}\\
& = & (3/4) \sup_{\theta_j = \pm 1}
     \Bigl\|\sum_j \theta_j \widetilde{R}_j \Bigr\|
\ge (1/4) \sum_j \| \widetilde{R}_j \|\\
& \ge & (1/8) \sum_j \| {R}_j \|
\ge  (1/8) \|\sum_j  {R}_j \| = (1/8) \|I_{Z_k}\| = 1/8,
\end{eqnarray*}
completing the proof.
\qed

Finally let us introduce  notations connected
with partitions of the set of natural numbers $\NN$,
which are essential in the sequel.
A subset $A \subset \NN$ is called an
interval if it is of the form
$A = \{i\,\mid k \le i \le n\}$.
Sets $A_1$ and $A_2$ are called consecutive intervals
if $\max A_i < \min A_j$, for $i, j = 1, 2$ and $i \ne j$.
A family of mutually disjoint subsets
$\Delta = \{A_{ m}\}_m$ is a partition of $\NN$,
if $\bigcup_m A_m  = \NN$.

For a partition  $\Delta = \{A_{ m}\}_m$ of $\NN$,
by ${\cal L}(\Delta)$ we denote the family
\begin{equation}
\label{suc_part_c}
{\cal L}(\Delta) = \{ L \subset \NN\, \mid\,
|L \cap A_m|=1 \ \mbox{for\ } m=1, 2, \ldots \}.
\end{equation}
If  $\Delta' = \{A_{ m}'\}_m$ is
another partition of $\NN$, we say that
$\Delta \succ \Delta'$,
if there exists a partition
${\cal J} (\Delta', \Delta)=  \{J_{ m}\}_m$ of $\NN$
such that
\begin{equation}
\label{suc_part_a}
\min J_m < \min J_{m+1}
\qquad \mbox{and} \qquad
A_{ m}' = \bigcup_{j \in J_m} A_{ j}
\qquad \mbox{for\ } m=1, 2,\ldots.
\end{equation}
In such a situation, for $m=1, 2, \ldots$,
${\cal K}(A_m', \Delta)$  denotes the family
\begin{equation}
\label{suc_part_b}
{\cal K}(A_m', \Delta) = \{ K \subset A_{ m}'\, \mid\,
|K \cap  A_{ j}|=1 \ \mbox{for\ } j \in J_m\}.
\end{equation}

Finally, if  $\Delta_i = \{A_{i, m}\}_m$,
for $i= 1, 2,\ldots$,  is a sequence of partitions of $\NN$,
with $\Delta_{1} \succ \ldots \succ \Delta_{i} \succ \ldots$,
we set, for $m=1, 2, \ldots$ and $i= 2,3, \ldots$
\begin{equation}
{\cal K}_{i, m}= {\cal K} (A_{i, m}, \Delta_{i-1}).
  \label{k_im}
\end{equation}

\smallskip
\section{General construction of subspaces
without local unconditional structure}

We will now present an abstract setting
in which it is possible to construct spaces
without local unconditional structure, but which still admit
a Schauder basis.
As it is quite natural, we work inside
a direct sum of several Banach spaces
with bases, with each basis having  certain
unconditional property related to some
partitions of $\NN$. The construction
of a required subspace
relies on an interplay between a ``good'' behaviour
of a basis on members of the corresponding partition
and a ``bad'' behaviour on sets which select one point from
each member of the  partition. (Recall that the notation
${\cal K}_{i, m}$  used below
was introduced in (\ref{k_im}).)

\begin{thm}
   \label{finite_part}
Let $X= F_{1}\oplus\ldots\oplus F_{4}$
be a direct sum of  Banach spaces
of cotype $r$, for some $r < \infty$,
and let $\{f_{i,l}\}_{l}$ be a  normalized
monotone Schauder basis
in $F_i$, for $i=1, \ldots, 4$.
Let $\Delta_{1} \succ \ldots \succ \Delta_{4}$
be partitions of $\NN$,
$\Delta_i = \{A_{i, m}\}_m$  for $i=1, \ldots, 4$.
Assume that  there is $C \ge 1$ such that
for every $K \in {\cal K}_{i, m}$  with
$i=2,3, 4$ and  $m=1, 2, \ldots$,
 the  basis
$\{f_{s,l}\}_{l \in K}$ in $\restr{F_s}{K}$ is
$C$-unconditional, for $s=1,  \ldots,4$;
moreover,
there is $\widetilde{C} \ge 1$ such that
for  $i=1, 2, 3$  and  $m=1, 2, \ldots$
we have
\begin{equation}
\|I: \restr{F_i}{A_{i,m} }\to
     \restr{F_{i+1}}{A_{i, m}} \|\le \widetilde{C}.
  \label{4_2_small}
\end{equation}
Assume finally  that one of the following
conditions is satisfied:
\begin{description}
\item[(i)]
there is a sequence $0 < \delta_m <1$ with
$\delta_m \downarrow 0$  such that
for every $i= 1, 2, 3$ and  $m=1, 2, \ldots$
and  every  $K \in {\cal K}_{i+1, m}$ we have
\begin{equation}
\|I: \restr{D(F_{1}\oplus\ldots\oplus F_{i})}{K}
\to \restr{F_{i+1}}{K}\|  \ge \delta_m ^{-1};
\label{4_2_large}
\end{equation}
\item[(ii)]
there is a sequence $0 < \delta_m <1$ with
$\delta_m \downarrow 0$ and
$\sum_m \delta_m ^{1/2} = \gamma <\infty$
such that
for every $i= 1, 2, 3$ and  $m=1, 2, \ldots$
and  every  $K \in {\cal K}_{i+1, m}$ we have
\begin{equation}
\|I: \restr{ F_{i+1}}{K}
    \to \restr{F_{i}}{K}\|  \ge \delta_m ^{-1}.
\label{4_3_large}
\end{equation}
\end{description}
Then there exists a subspace $Y$ of  $X$
without local  unconditional structure,
but which still admits a Schauder basis.
\end{thm}

\noindent{\bf Remarks\ }
{\bf 1.\,}
The space $Y$ will be constructed to have
a 2-dimensional Schauder decomposition.
If the bases $\{f_{i,l}\}_{l}$ are unconditional,
for $i=1, \ldots, 4$,
this decomposition will be  unconditional.

{\bf 2.\,}
Recall that a space which admits a
$k$-dimensional unconditional decomposition
has the GL-property (\cf\ [J-L-S])
(with the GL-constant depending on $k$).
Therefore the subspace $Y$ discussed in Remark 1
above has the GL-property but fails having
the local unconditional structure.

\medskip
\proof
We will define 2-dimensional subspaces
$Z_k$ of $X$ which will form a Schauder decomposition of
$Y = \overline{\spn}\,   [Z_k]_{k}$.
This decomposition will be $C'$-unconditional
on subsets associated with the partitions
$\Delta_1, \ldots, \Delta_4$,
for some $C'$ depending on $C$.
We shall use Proposition~\ref{borz}
to show that if $Y$  had the
local unconditional structure then, letting
$\psi = \psi(r, C_r(X))$ to be the function
defined in this proposition,
we would have
\begin{equation}
 \lust (Y) \ge \kappa \delta_t^{-\alpha}
  \label{contra_1}
\end{equation}
for an arbitrary $t=1, 2, \ldots$;
in case (i) we have
$\kappa > ( 2^{14} 3^3 C^4 \widetilde{C}^2 \psi  ) ^{-1}$
and  $\alpha = 1/3$;
in case (ii) we have
$\kappa > ( 2^{13} (1 + 4 \ga) C^3 \widetilde{C}^2 \psi  ) ^{-1}$
and $\alpha =1/2$.
This is impossible, which
will conclude the proof.

For $k=1, 2, \ldots$,  vectors $x_k$ and $y_k$
spanning $Z_k$ will be of the form
\begin{eqnarray*}
  x_k &=& \al_{1, k} f_{1, k}+ \ldots + \al_{4, k} f_{4, k},\\
  y_k &=& \al_{1, k}' f_{1, k}+ \ldots + \al_{4, k}' f_{4, k},
\end{eqnarray*}
such that  for $k=1, 2, \ldots$ and any scalars $s$ and $t$,
we will have
\begin{equation}
(1/2)\max (|s|, |t|) \le \|s x_k + t y_k \| \le 4 (|s| +|t|).
  \label{z_k_basis}
\end{equation}

Set $Z_k = [x_k, y_k]$, for $k=1, 2, \ldots$.
Clearly, $\{Z_k\}_k$ is a 2-dimensional
Schauder decomposition for $Y$,
in particular  $Y$ has a basis. Moreover,
for every  $i=2, 3, 4$ and
$m=1, 2, \ldots$ and   every $K \in {\cal K}_{i, m}$,
the decomposition  $\{Z_{k}\}_{ k \in K}$
is  $4 C$-unconditional.

Assume that $Y$ has the local unconditional structure.
Let $T$ be an operator obtained in
Proposition~\ref{borz}.
In particular, $T$ satisfies
(\ref{borz_ii})
for every $K \in {\cal K}_{i, m}$,
 and every  $i=2,3, 4$ and
$m=1, 2, \ldots$..
Let
$$
\left( \begin{array}{cc}
a_k & b_k\\
c_k & d_k
\end{array} \right)
$$
denote the matrix of $\restr{T}{{Z_k}}$ in the
basis $\{x_k, y_k\}$,
for $k=1, 2, \ldots$,\ie\ we have
$T(sx_k + t y_k) = (s a_k + t b_k ) x_k +
(s c_k + t d_k ) y_k$.
Comparing the operator norm of a $2 \times 2$
matrix with the $l_\infty^4$-norm of the
sequence of entries, and using
(\ref{z_k_basis}), we get that condition (iii)
of  Proposition~\ref{borz} implies that, for all
$k=1, 2, \ldots$,
\begin{equation}
\inf_\la \max (|a_k - \la|, |d_k - \la|, |b_k|, |c_k|)
\ge 2^{-5} \inf_\la \|\restr{T}{{Z_k}} - \la I_{Z_k}\|
\ge 2^{-8}.
  \label{far_id}
\end{equation}

\smallskip
For the rest of the argument we consider cases (i) and (ii)
separately. We start with (i).
Let $\ga_m = \de_m^{1/3}$, for $m=1, 2, \ldots$.
For $k \in A_{4,t}$, with $t=1, 2, \ldots$, put
\begin{equation}
\renewcommand\arraystretch{1.3}
\begin{array}{cccccc}
x_k & = & f_{1, k}& & + \ga_t f_{3, k}&  +  \ga_t^2 f_{4, k}\\
y_k & = & &  f_{2, k} & &  +  \ga_t^2 f_{4, k}.
\end{array}
  \label{def_4_2}
\end{equation}

Obviously, (\ref{z_k_basis}) is satisfied.
Fix an arbitrary $t=1, 2, \ldots$. For $i=1, 2, 3$,
let ${\cal M}_i= \{m\,|\, A_{i, m} \subset A_{4, t}\}$.
Note that (\ref{suc_part_a}) yields that
$\min  {\cal M}_i \ge t$ for  $i=1, 2, 3$.

For every $m \in {\cal M}_2$ pick $B \in {\cal K}_{2, m}$.
By (\ref{4_2_large}) we have
$$
\|I: \restr{F_1}{B} \to \restr{F_2}{B} \| \ge \ga_m^{-3};
$$
on the other hand, $\|f_{1, l}\|= \|f_{2, l}\|= \|I(f_{1, l})\|$.
By continuity, there exists a sequence $\{\bt_k\}_{k \in B}$
such that
$\|\sum_{k \in B} \bt_k f_{1, k}\|=1$
and
$\|\sum_{k \in B} \bt_k f_{2, k}\|= \ga_t^{-1}$.
Then, by (\ref{4_2_small}) and (\ref{def_4_2}) we have
\begin{eqnarray*}
\Bigl\|\sum_{k \in B} \bt_k x_k\Bigr\|
   &\le & \Bigl\|\sum_{k \in B} \bt_k f_{1, k}\Bigr\|
  +  \ga_t \Bigl\|\sum_{k \in B} \bt_k f_{3, k}\Bigr\|
  + \ga_t^2 \Bigl\|\sum_{k \in B} \bt_k f_{4, k}\Bigr\| \\
   &\le& 1 + (\ga_t \widetilde{C} + \ga_t^2 \widetilde{C}^2)
       \Bigl\| \sum_{k \in B} \bt_k f_{2, k}\Bigr\|
          \le 3 \widetilde{C}^2,
\end{eqnarray*}
while
\begin{eqnarray*}
\lefteqn{\Bigl\| T ( \sum_{k \in B} \bt_k x_k)\Bigr\|
   = \Bigl\|\sum_{k \in B} \bt_k (a_k x_k+ c_k y_k)\Bigr\|}\\
  &\ge& \Bigl\|\sum_{k \in B} \bt_k c_k f_{2,k}\Bigr\|
    \ge C^{-1} \inf_{k \in B}\,|c_k|\,
           \Bigl\|\sum_{k \in B} \bt_k f_{2,k}\Bigr\|
    \ge  C^{-1} \ga_t^{-1}  \inf_{k \in A_{2,m}}\,|c_k|.
\end{eqnarray*}

This implies, by (\ref{borz_ii}),
that for every $m \in {\cal M}_2$ there exists
$l \in A_{2, m}$ such that
$|c_l|\le  3 \, 4^2  C^3 \widetilde{C}^2  \psi \, \ga_t \lust (Y)$.
Denote the set
of these $l$'s by $L_2$ and observe
that $L_2 \in \restr{{\cal L}(\De_2)}{{\cal M}_2}$.
If we had $|c_l| > 2^{-10}$ for some $l \in L_2$,
then (\ref{contra_1}) would follow.
Therefore assume that
$|c_l| \le 2^{-10}$ for all $l \in L_2$.

For every $m \in {\cal M}_3$, set $B = L_2 \cap A_{3, m}$.
Then $B \in {\cal K}_{3, m}$
and by (\ref{4_2_large}) there exists
a sequence $\{\bt_k\}_{k \in B}$
such that
$$
\biggl\|\sum_{k \in B} \bt_k \biggl(\frac{f_{1, k}+ f_{2, k}}
  {\|f_{1, k}+ f_{2, k}\|}\biggr)\biggr\|=1
\qquad \mbox{and} \qquad
\Bigl\|\sum_{k \in B} \bt_k f_{3, k}\Bigr\|= \ga_t^{-2}.
$$
Observe that the basis
$\{( f_{1, k}+ f_{2, k})/ \|f_{1, k}+ f_{2, k}\|\}_{k \in B}$
is $2 C$-unconditional for every $B \in {\cal K}_{3, m}$.
Thus,
$$
\Bigl\|\sum_{k \in B} \bt_k f_{2, k}\Bigr\|
\le \Bigl\|\sum_{k \in B} \bt_k (f_{1, k} + f_{2, k})  \Bigr\|
 \le  4C.
$$
Hence,
$$
\Bigl\|\sum_{k \in B} \bt_k y_k\Bigr\|
   \le  \Bigl\|\sum_{k \in B} \bt_k f_{2, k}\Bigr\|
  + \ga_t^{2} \Bigl\|\sum_{k \in B} \bt_k f_{4, k}\Bigr\|
      \le  4\,C + \widetilde{C},
$$
and
\begin{eqnarray*}
\Bigl\| T \Bigl( \sum_{k \in B} \bt_k y_k\Bigr)\Bigr\|
   &=& \Bigl\|\sum_{k \in B} \bt_k (b_k x_k+ d_k y_k)\Bigr\| \\
  &\ge& \ga_t \Bigl\|\sum_{k \in B} \bt_k b_k f_{3,k}\Bigr\|
    \ge C^{-1} \ga_t^{-1} \inf_{k \in L_2 \cap A_{3,m}}\,|b_k|.
\end{eqnarray*}

Therefore, using  (\ref{borz_ii}) again,
for every $m \in {\cal M}_3$  pick
$l \in L_2 \cap A_{3, m}$ such that
$|b_l|\le 4^2 (4 C + \widetilde{C})   C^3 \psi \ga_t  \lust (Y)$.
Denote the set
of these $l$'s by $L_3$ and assume
as before that
$|b_l| \le 2^{-10}$ for all $l \in L_3$.
Moreover, $L_3 \subset L_2$
and $L_3 \in \restr{{\cal L}(\De_3)}{{\cal M}_3}$.

Finally, consider $K = L_3 \cap A_{4, t} \in {\cal K}_{4, t}$
and pick a sequence $\{\bt_k\}_{k \in K}$
such that
$$
\biggl\|\sum_{k \in K} \bt_k \biggl(\frac{f_{1, k}+ f_{2, k}+ f_{3, k}}
  {\|f_{1, k}+ f_{2, k}+ f_{3, k}\|}\biggr)\biggr\|=1
\qquad \mbox{and} \qquad
\Bigl\|\sum_{k \in K} \bt_k f_{4, k}\Bigr\|= \ga_t^{-3}.
$$
Since  $\{(f_{1, k}+ f_{2, k}+ f_{3, k})/
     \|f_{1, k}+ f_{2, k}+ f_{3, k}\|\}_{k \in K}$
is $3 C$-unconditional, for every $K \in {\cal K}_{4, t}$,
we have, for $i=1, 2, 3$,
$$
 \Bigl\| \sum_{k \in K} \bt_k f_{i, k} \Bigr\|
\le \Bigl\| \sum_{k \in K} \bt_k (f_{1, k}+ f_{2, k} + f_{3, k} )\Bigr\|
\le 3^2 C.
$$
Thus,
$$
\Bigl\|\sum_{k \in K} \bt_k (x_k - y_k)\Bigr\|
    \le \Bigl\| \sum_{k \in K} \bt_k f_{1, k} \Bigr\|
      + \Bigl\| \sum_{k \in K} \bt_k f_{2, k} \Bigr\|
          +  \ga_t \Bigl\| \sum_{k \in K} \bt_k f_{3, k} \Bigr\|
     \le 3^3 \,C.
$$
Moreover, since $|c_k|\le 2^{-10}$ and $|b_k|\le 2^{-10}$ for
$k \in L_3$, by (\ref{far_id}) we have
\begin{equation}
|a_k - b_k + c_k - d_k| \ge 2^{-9}
\qquad \mbox{for} \quad k \in L_3.
\label{paris2}
\end{equation}
Therefore
\begin{eqnarray*}
  \Bigl\| T \Bigl( \sum_{k \in K} \bt_k (x_k -y_k)\Bigr)\Bigr\|
&=& \Bigl\|\sum_{k \in K} \bt_k \Bigl((a_k - b_k) x_k
        +(c_k- d_k) y_k\Bigr)\Bigr\| \\
&\ge&  \ga_t^2 \Bigl\|\sum_{k \in K}
    \bt_k ((a_k - b_k) +(c_k- d_k) ) f_{4, k}\Bigr\|\\
&\ge&   C^{-1} 2^{-9} \ga_t^{-1}.
\end{eqnarray*}
Using  (\ref{borz_ii})  once more we get
$  3^3 4^2 C^3 \psi \, \lust (Y) \ge  C^{-1} 2^{-9} \ga_t^{-1}$,
which implies  (\ref{contra_1}).
This completes the proof of case (i).

\smallskip
In case (ii) the proof is very similar
and let us describe   necessary modifications.
Set $\ga_m = \de_m ^{1/2}$ for $m=1, 2, \ldots$.
For $k=1, 2, \ldots$ and
$k \in A_{2, m} \cap A_{3, s} $,
for some $m=1, 2, \ldots$ and  $s=1, 2, \ldots$,
set
\begin{equation}
\renewcommand\arraystretch{1.3}
\begin{array}{cccccc}
x_k & = & & \ga_{s} f_{2, k} & + f_{3, k}&  + f_{4, k}\\
y_k & = & \ga_{m} f_{1, k} &  & + f_{3, k}. &
\end{array}
  \label{def_4_3}
\end{equation}

Again, (\ref{z_k_basis}) is satisfied.
Fix an arbitrary  $t=1, 2, \ldots$, and define
${\cal M}_i$, for $i=1, 2, 3$ as before.
Using the fact that
$\|I: \restr{F_2}{K} \to   \restr{F_1}{K} \|\ge \ga_m^{-2}$,
for every $K \in {\cal K}_{2, m}$ and every $m \in {\cal M}_2$,
one can show, using (\ref{4_2_small}) and
(\ref{borz_ii}) in a similar way as before,
that there is a set
$L_2 = \{l_m\}_{m \in {\cal M}_2}
\in \restr{{\cal L}(\De_2)}{{\cal M}_2}$
such that
\begin{equation}
|c_{l_m}| \le 3 \, 4^2 C^3 \widetilde{C}^2 \psi \ga_m    \lust (Y)
\qquad \mbox{for\ }m \in {\cal M}_2.
  \label{4_3_c}
\end{equation}
One can additionally assume that
$|c_{l_m}| \le 2^{-10}$ for all
$m \in {\cal M}_2$, otherwise,
since $\min {\cal M}_2 \ge t$ implies
$ \ga_m \le \ga_t $, we would immediately
get (\ref{contra_1}) with $\alpha = 1/2$..

Now for every $s \in {\cal M}_3$ consider the set
$B= L_2 \cap A_{3, s} \in {\cal K}_{3, s}$, and pick
a sequence $\{\bt_k\}_{ k \in B}$
such that
$$
\Bigl\|\sum_{ k \in B} \bt_k f_{3, k}\Bigr\|=1
\qquad \mbox{and} \qquad
\Bigl\|\sum_{ k \in B} \bt_k f_{2, k}\Bigr\|\ge \ga_s^{-2}.
$$
If ${\cal M}_{2, s}$ denotes the set of indices
$m\in {\cal M}_2$ such that
$l_m \in L_2 \cap A_{3, s}= B$, then
\begin{equation}
\label{paris1}
\Bigl\|\sum_{ k \in B} \bt_k y_{ k}\Bigr\|
\le \sum_{m \in {\cal M}_{2, s}} \ga_m |\bt_{l_m}| +1 \le 2\ga+1,
\end{equation}
where the first term in the estimate
is obtained  by first
using the triangle inequality and then using the fact
that  since
$\{f_{3, l_m}\}_{m \in {\cal M}_{2, s}}$
is a monotone  basic sequence, then
$ |\bt_{l_m}|  \le 2 $ for all $l_m \in B$.

We also have
\begin{eqnarray*}
\Bigl\| T \Bigl( \sum_{k \in B}
               \bt_k y_k\Bigr)\Bigr\|
  &\ge& \ga_{s}\Bigl\|\sum_{k \in B}
                \bt_k b_k f_{2,k}\Bigr\|\\
   & \ge& C^{-1} \ga_{s} \inf_{k \in B}\,|b_k|\,
        \Bigl\|\sum_{ k \in L_2 \cap A_{3, s}} \bt_k f_{2, k}\Bigr\|\\
  & \ge&  C^{-1} \ga_s^{-1} \inf_{k \in L_2 \cap A_{3,s}}\,|b_k|.
\end{eqnarray*}
Thus there exists a set $L_3 \in \restr{{\cal L}(\De_3)}{{\cal M}_3}$,
$L_3 = \{l_s'\}_{s \in {\cal M}_3}$, such that  $L_3 \subset L_2$
and
\begin{equation}
|b_{l_s'}| \le 4^2 ( 2\ga +1) C^3 \psi \ga_s   \lust (Y)
\qquad \mbox{for\ }s \in {\cal M}_3;
  \label{4_3_b}
\end{equation}
and since $\min {\cal M}_3 \ge t$ implies
$\ga_s \le \ga_t$,
one can additionally assume  that
$|b_{l_s'}| \le 2^{-10}$, for all
$s \in {\cal M}_3$.

Finally set
 $K = L_3 \cap A_{4, t} \in {\cal K}_{4, t} $.
Pick $\{\bt_k\}_{k \in K}$ such that
$\|\sum_{ k \in K} \bt_k f_{4, k}\|=1$
and
$\|\sum_{ k \in K} \bt_k f_{3, k}\| \ge \ga_t ^{-2}$.
Then, by the triangle inequality and by the monotonicity
of the basis $\{f_{4,k} \}_k$ we get,
similarly as in (\ref{paris1}),
$$
\Bigl\|\sum_{k \in K} \bt_k (x_k - y_k)\Bigr\|
   \le  2\sum_{m \in {\cal M}_2} \ga_m
     +  2\sum_{s \in {\cal M}_3} \ga_s
     + \Bigl\| \sum_{k \in K} \bt_k  f_{4, k}\Bigr\|
        \le 1 + 4 \ga.
$$
On the other hand,
 by (\ref{far_id}),  (\ref{4_3_c}) and (\ref{4_3_b})
we again have (\ref{paris2}).
Thus
\begin{eqnarray*}
  \Bigl\| {T} \Bigl( \sum_{k \in K} \bt_k (x_k -y_k)\Bigr)\Bigr\|
&\ge& \Bigl\|\sum_{k \in K}
    \bt_k (a_k - b_k ) + (c_k - d_k)  f_{3, k}\Bigr\|\\
&\ge& C^{-1} 2^{-9} \Bigl\|\sum_{k \in K} \bt_k f_{3, k}\Bigr\|
\ge   C^{-1} 2^{-9} \ga_t^{-2}.
\end{eqnarray*}

 Using (\ref{borz_ii}) we get
 $ \lust (Y) \ge
 (2^{13}(1 + 4 \ga ) C^3  \psi ) ^{-1}\ga_t^{-2}$,
hence  (\ref{contra_1}) follows,
completing the proof of case (ii).
\qed



\smallskip
\section{Subspaces of spaces with unconditional basis}

Our main application of the construction of
Theorem~\ref{finite_part} is the
following result on
subspaces of spaces with unconditional basis.

\begin{thm}
  \label{l_2_sub}
Let $X$ be a Banach space with an unconditional basis
and of cotype $r$, for some $r < \infty$.
If $X$ does not contain a subspace isomorphic  to $l_2$ then
there exists a subspace $Y$ of  $X$ without local unconditional
structure, which admits a Schauder basis.
\end{thm}

In particular, every Banach space of cotype $r$,
for some $r < \infty$, contains either $l_2$ or a
subspace without unconditional basis.

We present  now  the proof of the theorem, leaving
corollaries and further applications to the  next
section.

The  argument
is based on a construction, for a
given Banach space  $X$,
of a    direct sum  inside $X$
of  subspaces $F_i$ of $X$,
and of   partitions $\Delta_i$ of $\NN$
such that
Theorem~\ref{finite_part} can be applied.
This construction requires several steps.

The first  lemma is  a simple generalization
to finite-dimensional
lattices of the fact that the  Rademacher functions
in $L_p$ are equivalent to the
standard unit vector basis in $l_2$.


\begin{lemma}
\label{orbit}
Let $E$ be an $N$-dimensional Banach space
with a 1-uncon\-di\-tio\-nal  basis $\{e_{j}\}_j$
and for $2 \le r < \infty$ let
$C_r(E)$ denote the cotype $r$ constant
of $E$.
If $ m \le \log_2 N$ then there exist
normalized vectors $f_1, \ldots, f_m$
in $E$, of the form
\begin{equation}
\label{orbit_a}
f_l = \sum_j \ep_j^{(l)} \al_j e_j
\qquad \mbox{ for\ }
l=1, \ldots, m,
\end{equation}
for some sequence  of scalars $\{\al_j\}$
and $\ep_j^{(l)}= \pm 1$ for
$l=1, \ldots, m$ and $j=1, \ldots, N$;
and such that
\begin{equation}
\label{orbit_b}
de\Bigl(\spn [f_l],\, l_2^m\Bigr) \le C,
\end{equation}
where $C$  depends on $r$ and on
the cotype $r$ constant of $E$.
\end{lemma}
\proof
Since $E$ is a discrete Banach lattice, the cotype $r$
assumption implies that $E$ is $q$-concave, for any $q > r$
(\cf\  [L-T.2]).
Setting \eg\ $q = 2r$, the $q$-concavity constant of $E$
depends on $r$ and  $C_r(E)$. By a lattice renorming
we may and will
assume that this constant  is equal to 1
(\cf\  [L-T.2] 1.d.8);
the general case will follow by  adjusting $C$.

For $1 \le p < \infty$, let $\|\cdot\|_{L_p}$
be the norm defined on $\Rn{N}$ by
$ \|t\|_{L_p} =  (N^{-1} \sum_{j=1}^N |t_j|^p )^{1/p}$,
for  $t= (t_j)\in \Rn{N}$.
It is well known consequence of Lozanovski's theorem
(see [T], 39.2 and 39.3 for a related result)
that there exist $\al_j >0$, $j=1, \ldots, N$,
such that
\begin{equation}
\label{lozanov}
\|t\|_{L_1}
\le \|  \sum_{j=1}^N \al_j t_j e_j \|
\le  \|t\|_{L_q}
\qquad \mbox{for\ }
t= (t_j)\in \Rn{N}.
\end{equation}
%

Fix an integer  $m \le \log_2 N$.
By Khintchine's inequality  there exist
vectors $r_l = \{r_l(j)\}_{j=1}^N$, with
$r_l(j)= \pm 1$ for $j=1, \ldots, N$, $l=1, \ldots, m$,
such that  for every $(b_l)\in \Rn{m}$ we have
\begin{equation}
\label{kh}
2^{-1/2} (\sum_{l=1}^m |b_l|^2 )^{1/2}
\le \| \sum_{l=1}^m b_l r_l\|_{L_1}
\le \| \sum_{l=1}^m b_l r_l\|_{L_q}
\le C_q (\sum_{l=1}^m |b_l|^2 )^{1/2}.
\end{equation}
Setting $f_l = \sum_{j=1}^N  r_l(j) \al_j e_j$, for $l=1, \ldots, m$,
we get, by (\ref{lozanov}),
$$
\| \sum_{l=1}^m b_l r_l\|_{L_1}
\le   \| \sum_{l=1}^m b_l f_l\|
= \|\sum_{j=1}^N \al_j (\sum_{l=1}^m b_l r_l(j))  e_j\|
\le  \| \sum_{l=1}^m b_l r_l\|_{L_q},
$$
for every $(b_l)\in \Rn{m}$.
This combined with (\ref{kh}) completes the required
estimate.
\qed

\rem
As it was pointed out to us by B.~Maurey, Lemma~\ref{orbit}
could be replaced by the contruction of L.~Tzafriri [Tz],
which implies the existence of a function $\varphi(N)$, with
$\varphi(N) \to \infty$ as $N \to \infty$, such that for
$m \le \varphi(N)$
every $N$-dimensional space $E$ as in the lemma
contains   normalized vectors $f_1, \ldots, f_m$
satisfying  (\ref{orbit_b}),
which  are of the form
$f_l = \alpha \sum_j \pm e_j$, with an appropriate constant $\alpha$.

\medskip

The next proposition is the key for our argument.  To simplify the
statement, let us introduce one more notation.  Given a
partition $\Delta= \{A_{ m}\}_m$ of $\NN$
into consecutive intervals
and a  space  $F$  with  a normalized Schau\-der basis
$\{f_{l}\}_{l}$ and $C \ge 1$, we call a pair
$\{\Delta, F\}$
$C$-regular, if the following conditions
are satisfied:
\begin{description}
\item[(i)] $ de \Bigl(\restr{F }{A_{m} },\, l_2 ^{|A_{ m} |}\Bigr)
  \le   C $
for $m=1, 2, \ldots$;
\item[(ii)]
  for every $L \in {\cal L}(\De) $, the basis
  $\{f_{l}\}_{l \in L}$ in $\restr{F }{L}$ is
  1-unconditional (here $ {\cal L}(\Delta )$ is as in
  (\ref{suc_part_c}));
\item[(iii)] for arbitrary $L, L' \in {\cal L}(\Delta ) $ one has
$de \Bigl(\restr{F }{L}, \restr{F }{L'}\Bigr) =1$.
\end{description}

Observe that condition (iii) means that if $L= \{l_m\}_m$, $L'=
\{l_m'\}_m$, with $l_m, l_m' \in A_{m} $ for $m=1, 2,\ldots$, then
for every sequence of scalars $(b_m)$ one has
\begin{equation}
  \Bigl\|\sum_{m} b_m f_{l_m}\Bigr\| = \Bigl\|\sum_{m} b_m
  f_{l_m'}\Bigr\|.    \label{three}
\end{equation}

\begin{prop}
\label{seq_partitions}
Let $E_1, E_2 \ldots $ be Banach spaces
of cotype $r$, for some $r < \infty$.
Let  $\{e_{i,j}\}_j$ be a
1-un\-con\-di\-tio\-nal ba\-sis in $E_i$,
and assume that
no  sequence of disjointly supported vectors
in  $E_1\oplus\ldots\oplus E_i$
is equivalent to the standard unit vector basis in $l_2$,
for $i=1, 2,\ldots$.
Then there exists $C$, depending on $r$
and the cotype $r$ constants of $E_i$,
such that
there exist  subspaces $F_i \subset E_i$
with  normalized Schauder bases
$\{f_{i,l}\}_{l}$,
and partitions  $\Delta_{i} = \{A_{i, m}\}_m$
of $\NN$ into consecutive intervals,
for $i=1, 2, \ldots$,
with $\Delta_{1} \succ \Delta_{2} \succ \ldots$,
satisfying the following:
for each $i=1, 2, \ldots$
$\{\Delta_i, F_i\}$  is $C$-regular and
one of the following mutually exclusive conditions
is satisfied:
%
either for every ${L \in {\cal L}(\Delta_i)}$ one has
\begin{equation}
 \|I: l_2 \to \restr{F_i}{L}\|  =\infty,
  \label{one}
\end{equation}
or  for every ${L \in {\cal L}(\Delta_i)}$ one has
\begin{equation}
 \|I: l_2 \to \restr{F_i}{L}\|  <\infty.
  \label{two}
\end{equation}
%
Furthermore, one also has
\begin{description}
\item[(iv)]
If  (\ref{one}) holds for some $i$, then the
partition $\Delta_{i+1} = \{A_{i+1, m}\}_m$
satisfies
\begin{equation}
{\inf_m} \,  \inf\, \Bigl\{
   2^{-3m}\|I: l_2^{|K|} \to \restr{F_i }{K} \,\|\,
     \Bigm|  K \in {\cal K}_{i+1, m} \Bigr\} \ge C.
  \label{one_finite}
\end{equation}
On the other hand, let  $M$ denote  the set
(which may be empty)
of all $ s \in \NN$  such that
for every ${L \in {\cal L}(\Delta_s)}$ one has
$  \|I: l_2 \to \restr{F_s}{L}\|  <\infty $.
If $i \in M$, put $M_i = M \cap \{1, \ldots, i\}$;
then  the partition $\Delta_{i+1} = \{A_{i+1, m}\}_m$
satisfies
\begin{eqnarray}
&&\inf_m \,  \inf\, \Bigl\{
 2^{-3m} \|I:
     \restr{D\Bigl( \sum_{s \in M_i}\oplus F_s\Bigr)}{K}
    \to l_2 ^{|K|}\|  \, \Bigm|\nonumber\\
  &&\qquad \qquad \qquad \qquad
            \qquad K \in {\cal K}_{i+1, m} \Bigr\} \ge C.
  \label{two_finite}
\end{eqnarray}
\end{description}
\end{prop}
\proof
In the first part of the proof we
show that given space $E$ of cotype $r$
with a
1-un\-con\-di\-tio\-nal  ba\-sis
$\{e_j\}_j$,
and a partition $\Delta =  \{A_{ m}\}_m$   of $\NN$
into consecutive intervals,
there exists a subspace
$F \subset E$ with a normalized Schauder basis
$\{f_{l}\}_{l}$ such that $\{\Delta, F\}$
is $C$-regular, for an appropriate
constant $C$, and that either (\ref{one})
or (\ref{two}) is satisfied for
every ${L \in {\cal L}(\Delta)}$.

For an arbitrary $m=1, 2,\ldots$,
let $k_{m}= |A_{ m}|$ and
let $E^{(m)} = \spn\{e_{ j}\,\mid\,
2^{k_{m}} < j \le 2^{k_{m+1}}\}$.
Since $\dim E ^{(m)}\ge 2^{k_{m}}$,
by  Lemma~\ref{orbit} there exist vectors
$f_{ l} \in E^{(m)} $, for
$l \in A_{ m}$, such that
\begin{equation}
\label{de_small}
de \Bigl(\spn [f_{ l}]_{ l \in A_{ m}},\, l_2^{k_m}\Bigr) \le C;
\end{equation}
and there is a sequence $\{\al_{ j}\}$ of real numbers
such that the $f_{ l}$'s are of the form
\begin{equation}
\label{form}
f_{ l}= \sum_{j = 2^{k_{m}}+1}^{2^{k_{m+1}}} \pm \al_{ j} e_{ j}
\qquad \mbox{for \ } l \in A_{ m},\, m=1, 2, \ldots.
\end{equation}

We let  $ F = \spn [f_{ l}]_{l}$.
Then (i) is implied by (\ref{de_small}).
Next observe that
$f_{ l}$ and $f_{ l'}$ have consecutive
supports, whenever
$l \in  A_{ m}$ and $l' \in  A_{ m'}$ and  $m \ne m'$.
This and (\ref{de_small}) easily yield that
$\{f_{ l}\}_l$  is a Schauder basis in $F$.
Also, $\{f_l\}_{l \in L}$ is a 1-unconditional basis
in $\restr{F}{L}$, for every $L \in {\cal L}(\De)$, which
shows (ii).

By (\ref{form}) we get that if $(b_m)$ is a scalar sequence
then for every $L = \{l_m\}_m \in {\cal L}(\De)$, the vector
$\sum_{m} b_m f_{l_m}$ is of the form
$$
\sum_{m} b_m \sum_{j = 2^{k_{m}}+1}^{2^{k_{m+1}}} \pm \al_{ j} e_{ j};
$$
a specific choice of the $l_m$'s  which constitute the set $L$
effects only the choice of the signs in the  inner summation.
Since the basis $\{e_{j}\}$ is 1-unconditional,
(\ref{three}) follows, hence  (iii) follows as well.

Finally observe that for a  fixed $L \in {\cal L}(\De) $,
exactly one of conditions
(\ref{one}) and  (\ref{two}) holds.
Moreover, by (iii), the norms of the formal identity operators
involved do not depend on a choice of the set
$L \in {\cal L}(\De) $.

\medskip
We now pass to the second part of the proof, the inductive
construction of $\Delta_i$'s and $F_i$'s, which
ensures also condition (iv).
Let $A_{1, m}= \{m\}$ for $m=1, 2, \ldots$ and let
$\Delta_1= \{A_{1, m}\}_m$.

Assume that $i \ge 1$ and that partitions
$\Delta_{1} \succ \ldots \succ \Delta_{i}$
and subspaces $F_1, \ldots, F_{i-1}$ have
been constructed
to satisfy conditions (i)--(iv).
Let $F_i \subset E_i$ be a subspace
constructed in the first part of the proof
for $\Delta = \Delta_i$.
The construction of $ \Delta_{i+1}$ depends
on which of two,
(\ref{one}) or  (\ref{two}), holds for $F_i$.

Assume first that (\ref{one}) holds
and fix an arbitrary set $L  \in {\cal L}(\De_i)$.
Enumerate $L = \{l_j\}_j$ with $l_j \in A_{i, j}$
for $j=1, 2, \ldots$. There exist
$1= j_0 < j_1 < \ldots < j_m < \ldots$
such that if $J_m = \{j_{m-1}\le j < j_m\}$,
then
\begin{equation}
 \|I: l_2 ^{|J_m|} \to \restr{F_i}{\restr{L}{J_m}}\| \ge  C^2 2^{3m}
\qquad \mbox{for\ } m=1, 2, \ldots.
  \label{one_constr}
\end{equation}
We then set
\begin{equation}
A_{i+1, m} = \bigcup_{j \in J_m} A_{i, j}
\qquad \mbox{for\ } m=1, 2, \ldots.
  \label{new_part}
\end{equation}
By (\ref{three}) and (\ref{one_constr}) it is clear
that (\ref{one_finite}) is satisfied in this case.

Assume now that (\ref{two}) holds, so $i \in M$.
There is  a constant $C'$ such that for all $s \in M_i$
the estimate $ \|I: l_2 \to \restr{F_s}{L}\|  < C'$
holds for all $L \in {\cal L}(\De_s)$; hence
also  for  all $L \in {\cal L}(\De_i)$,
since sets from ${\cal L}(\De_i)$ are subsets of sets
from ${\cal L}(\De_s)$, for every $s < i$.
Fix an arbitrary $L \in {\cal L}(\De_i)$.
We then have
$$
 \|I: l_2 \to
    \restr{D\Bigl(\sum_{s \in M_i}\oplus F_{s}\Bigr)}{L}\|
 < |M_i|\,C'.
$$
Note that if $l, l' \in L  \in {\cal L}(\De_i)$
and $l\ne l'$ then
$f_{s, l}$ and $f_{s, l'}$ have consecutive supports,
hence $\{f_{s,l}\}_{l \in L}$ forms a block basis of
$\{e_{s,j}\}_j$, for $s \in M_i$.
Therefore by our assumptions, the basis
$\{\sum_{s \in M_i} f_{s,l}\}_{l \in L}$
in $D(\sum_{s \in M_i}\oplus F_{s})$
is not equivalent to the standard unit
vector basis in $l_2$.
Thus
\begin{equation}
 \|I: \restr{D\Bigl(\sum_{s \in M_i}\oplus F_{s}\Bigr)}{L}
       \to  l_2\| = \infty.
  \label{inv_large}
\end{equation}

Now the construction of a partition $\Delta_{i+1}$ satisfying
(\ref{two_finite}) is done by formulas completely analogous to
(\ref{one_constr}) and (\ref{new_part}), in which the use of
(\ref{one}) is replaced by (\ref{inv_large}).
\qed

Finally, the proof of the main result follows formally
from Proposition~\ref{seq_partitions}.

\smallskip
\noindent{\bf Proof of Theorem~\ref{l_2_sub}\ }
Write $X$ as an unconditional sum
$X = \sum_i \oplus E_i$, of 13
spaces $E_i$,
each with a 1-unconditional basis $\{e_{i,j}\}_j$.
Let  $\Delta_{1} \succ \ldots \succ \Delta_{13}$
be partitions  of $\NN$
and $F_i  \subset E_i$ be subspaces
with  Schauder bases  $\{f_{i,l}\}_{l}$,
constructed  in Proposition~\ref{seq_partitions}.
Renorming the spaces $F_i$ if necessary,
we may assume that the bases  $\{f_{i,l}\}_{l}$
are monotone.

Now the $C$-regularity properties imply all the preliminary
assumptions of Theorem~\ref{finite_part}, including
(\ref{4_2_small}).
To prove  the remaining conditions (i) or (ii) observe that
either there exist four consecutive  spaces $\{F_{i_k}\}_k$ satisfying
(\ref{two}), or (\ref{one}) holds for some three
(not necessarily consecutive) spaces $\{F_{i_k}\}_k$.

In either  case, we  let  $\Lambda_k = \Delta_{i_k}$
and  $F_k' = F_{i_k}$, for $k=1, \ldots, 4$
(in the latter case we set $i_4 = i_3 +1$).

It is easy to check that in the former case,
(\ref{two}) yields (\ref{two_finite}), while in
the latter case (\ref{one}) yields
(\ref{one_finite}).
Thus the remaining assumptions of Theorem~\ref{finite_part}
are satisfied with  $\de_m = 2 ^{-3m}$,
which concludes the proof.
\qed

\smallskip
\section{Corollaries and further applications}

Recall a still open question  whether a Banach
space whose all subspaces have an unconditional basis
is isomorphic to a Hilbert space.
>From results on the  approximation property
by Enflo, Davie, Figiel and Szankowski,
combined with Maurey--Pisier--Krivine theorem, it follows
that  such a space $X$
has, for every $\ep >0$, cotype $2+\ep$
and type $2 - \ep$ (\cf\ \eg\ [L-T.2], 1.g.6).
Theorem~\ref{l_2_sub} obviously implies
that $X$ has  a  much stronger property:
its every  infinite-dimensional
subspace contains an isomorphic copy of $l_2$.
A space $X$ with this property is called
$l_2$-saturated.

\begin{thm}
\label{l_2_saturated}
Let  $X$ be an infinite-dimensional Banach space
whose all subspaces have an unconditional basis. Then
$X$ is   $l_2$-saturated.
\end{thm}


Another well known open problem, going back to
Mazur and Banach, concerns so-called  homogeneous
spaces.
An infinite-dimensional Banach space
is called homogeneous  if it is isomorphic to each of
its   infinite-dimensional subspaces.
The question is whether every homogeneous Banach space is
isomorphic to a Hilbert space.
The same general argument as before shows that
a homogeneous  space $X$
has   cotype $2+\ep$  and type $2 - \ep$,
for every $\ep >0$.
W.~B.~Johnson showed in
[J] that if both $X$ and $X^*$  are homogeneous
and $X$ has the Gordon--Lewis property, then
$X$ is isomorphic to a Hilbert space.
More information about homogeneous spaces the reader can find
in [C].
The following obvious corollary
removes the assumption on  $X^*$,
however  it requires  a stronger property
of $X$ itself.

\begin{thm}
  \label{homogen}
If a homogeneous Banach space $X$ contains an infinite
unconditional basic sequence then $X$ is isomorphic
to a Hilbert space.
\end{thm}

Let us recall here that it was believed for a long time
that every Banach space might  contain an infinite
unconditional basic sequence. This conjecture was disproved
only recently by W.~T.~Gowers and B.~Maurey in [G-M], who
actually constructed a whole class of Banach spaces
failing this and related properties.

\medskip

Let us now discuss some  easy consequences
of the main construction,
which might be of independent interest.

\begin{cor}
  \label{inf_sums}
Let $X= F_{1}\oplus\ldots\oplus F_{4}$
be a direct sum of  Banach spaces
of cotype $r$, for some $r < \infty$,
and assume that  $F_i$ has a 1-unconditional
basis  $\{f_{i,l}\}_{l}$.
for $i=1, \ldots, 4$.
Assume that the basis  $\{f_{i,l}\}_l$ dominates
$\{f_{i+1,l}\}_l$,
and  that no subsequence of  $\{f_{i,l}\}_l$ is equivalent
to the corresponding subsequence of  $\{f_{i+1,l}\}_l$,
for $i=1, 2, 3$.
Then there exists a subspace $Y$ of  $X$
without local unconditional structure, which admits
an unconditional decomposition into  2-dimensional
spaces.
\end{cor}
\proof
Let $\Delta_{1} \succ \ldots \succ \Delta_{4}$
be arbitrary partitions of $\NN$ into infinite
subsets $\{A_{i,m}\}$. The domination assumption
implies  (\ref{4_2_small}). On the other hand,
the second  assumption allows for
a construction of partitions which also
satisfy (\ref{4_3_large}).
Hence the conclusion
follows from Theorem~\ref{finite_part} and  Remark 2 above.
\qed

\rem
In fact, Corollary~\ref{inf_sums} can be  proved
directly from Proposition~\ref{borz}.
To define $x_k$ and $y_k$ spanning $Z_k$,
let $\La_2 = \{B_m\}_m$ be any partition
of $\NN$ into infinite sets and  write
each $B_m$  as a union $B_m = \bigcup_n B_{m, n}$
of an infinite number of infinite sets $ B_{m, n}$.
(Using the natural enumeration of $\NN\times \NN$, we
get this way a partition $\La_1 = \{B_{m, n}\}_{m, n}$
with $\La_{1} \succ \La_2$.) Then for $k \in B_{m, n}$,
with $m, n = 1, 2, \ldots$ put
$$
\renewcommand\arraystretch{1.3}
\begin{array}{cccccc}
x_k & = & & 2^{-m} e_{2, k} & + e_{3, k}&  + e_{4, k}\\
y_k & = & 2^{-m-n} e_{1, k} &  & + e_{3, k}. &
\end{array}
$$
The rest of the proof is the same as in
 case (ii) of Theorem~\ref{finite_part}.

\medskip

If $\{x_i\}$ is a basic sequence in a Banach space $X$, and
$1 \le p < \infty$, we say that $l_p$ is
crudely finitely sequence  representable  in $\{x_i\}$
if there is a constant $C \ge 1$ such that
for every $n$ there is a subset $B_n \subset \NN$ such that
$\{x_i\}_{i \in B_n}$ is $C$-equivalent to the unit vector basis
in $l_p^n$.

\begin{cor}
  \label{fin_rep_seq}
Let $X$ be a Banach space
of cotype $r$, for some $r < \infty$,
and with a  1-unconditional
basis  $\{e_{l}\}_{l}$;
let $1 \le p < \infty$.
Assume that no  sequence  $\{x_j\}_j$ of disjointly
supported vectors of the form
$x_j =  \sum_{l \in L_j} e_l$,
where $|L_j| \le 3$ for $j=1, 2,\ldots$,
is equivalent to the unit vector basis of $l_p$.
Moreover assume that
$X$ has one of the following properties:
\begin{description}
\item[(i)]
$l_p$ is crudely finitely sequence representable in
$\{e_l\}_l$, and the basis  $\{e_l\}_l$ either is dominated by
or dominates the standard unit vector basis in $l_p$;
\item[(ii)]
$l_p$ is crudely  finitely sequence
representable in every subsequence of $\{e_l\}_l$.
\end{description}
Then $X$ contains a subspace $Y$
without local unconditional structure, which admits
a 2-dimensional  unconditional decomposition.
\end{cor}
\proof
First observe a general fact concerning a basis $\{e_l\}_l$
whose no subsequence is
dominated by the standard unit vector basis in $l_p$.
An easy diagonal argument shows that if a partition
$\Delta = \{A_j\}_j$ of $\NN$ into finite
sets is given then for an arbitrary
$M$ and every $j_0\in \NN$
there is $j_1 > j_0$ such that for any set  $K \subset \NN$
such that $|K| =  j_1 - j_0$ and $|K \cap A_j| = 1$
for $ j_0 < j \le j_1$, one has
$ \|I:  l_p^{|K|} \to \restr{E}{K}\| \ge M$.
In particular, given constant $C$,
there exists a partition
$\Delta' = \{A'_m\}_m$ of $\NN$,
with $\De \succ \De'$ such that for every
$m=1, 2, \ldots$ and for every
$K \in {\cal K}(A'_m, \De)$ one has
\begin{equation}
\|I:  l_p^{|K|} \to \restr{E}{K}\| \ge C\,2^{3m}.
  \label{lower}
\end{equation}

Now,
in case (i),  write $X$ as a direct sum
$ E_{1}\oplus\ldots\oplus E_{4}$,
such that each $E_i$ has a 1-unconditional
basis  $\{e_{i,l}\}_{l}$.
Assume that the basis $\{e_l\}_l$ dominates
the basis in $l_p$, hence so does every
basis  $\{e_{i,l}\}_{l}$. Using the general
observation above,
we can define  by induction  partitions
$\Delta_1 \succ \ldots \succ \Delta_4$
and subsequences $\{f_{i,j}\}_{j}$ of $\{e_{i,l}\}_{l}$,
so that  for all $k$ and all
$A=A_{k,m}\in \Delta_k$,
sequences $ \{f_{i_k,j}\}_{j \in A }$
are $C$-equivalent to the standard unit
vector basis in  $l_p^{|A|}$,
and at the same time, the spaces
$\spn [f_{i_k,j}]_{j \in K}$,
with $K \in {\cal K}_{k+1, m}$,
satisfy the  lower estimate (\ref{lower}).
Thus   (\ref{4_3_large}) is satisfied
(with $\de _m = 2 ^{-3m}$).

If the basis $\{e_l\}_l$ is dominated
by the basis in $l_p$,  so is  every
basis  $\{e_{i,l}\}_{l}$, and also
all  bases in  $D(E_{1}\oplus\ldots\oplus E_{i})$,
for $i=1, 2, 3$.
An analogous  argument  as before, which additionally
requires the assumption on sequences $\{x_j\}$,
leads to a construction of partitions satisfying
(\ref{4_2_large}). Then
the existence of the subspace $Y$ follows
from Theorem~\ref{finite_part} and
Remark 1 after its statement.

In case (ii),   write
$X =  E_{1}\oplus\ldots\oplus E_{7}$.
By passing to subsequences  we get that
for each $i$, each subsequence of the  basis
$\{e_{i,l}\}_{l}$, either is dominated by  or dominates
the standard unit vector
basis in $l_p$, for $i=1, \ldots, 7$.
Therefore there is a set
$I= \{i_1,\ldots, i_4\}$ such that  for all $i \in I$,
the bases  $\{e_{i,l}\}_{l}$ have the same,
either former or latter, domination
property.  Then the proof can be concluded
the same way as in case (i).
\qed

For $ 1\le q < \infty$, the space
$L_q([0,1])$ contains a subspace isomorphic to
$X = (\sum_n\oplus l_2^n)_q$, which, for $q \ne 2$,
satisfies  the assumptions  of
Corollary~\ref{fin_rep_seq}~(i) for $p=2$.
Therefore $L_q([0,1])$ contains a subspace
without  local unconditional structure
but which admits a  2-dimensional
unconditional decomposition.
By Remark 2 in Section 2,
this subspace has the Gordon-Lewis  (GL-) property.
For $1 \le q <2$, this
gives a somewhat more elementary proof
of Ketonen's result [Ke]. For $ 2 < q < \infty$
the construction seems to be new.
Ketonen's result could be also derived from
Corollary~\ref{inf_sums} by noticing that
in this case the space $L_q([0,1])$
contains a subspace isometric to
$(l_{q_1}\oplus\ldots \oplus l_{q_4})_q$, for
$1 \le q\le q_1< \dots < q_4 < 2$
(\cf\ \eg\ [L-T.2], 2.f.5).

\medskip
Corollary~\ref{fin_rep_seq} can also
be applied to construct subspaces without
local unconditional structure in $p$-convexified
Tsirelson spaces $T_{(p)}$ and in their duals.
This solves the question left open
in [K].
The  spaces  $T_{(2)}$ and $T_{(2)}^*$
provide the most important
examples of so-called weak Hilbert spaces,
and they were discussed  in [P].
For general $p$ and notably for $p=1$,
these spaces were presented in detail
in  [C-S]. First construction
of a weak Hilbert space without unconditional basis
was given by R.~Komorowski in [K] by a method
preceeding the technique presented here.

\begin{cor}
  \label{tsir}
The $p$-convexified Tsirelson space $T_{(p)}$,
for $1 \le p < \infty$, and the
dual Tsirelson $T_{(p)}^*$,
for $1 < p < \infty$,
contain subspaces without
local unconditional structure, but which admit
2-di\-men\-sio\-nal unconditional
decomposition; in particular they have the Gordon--Lewis property.
\end{cor}
\proof
The spaces  $T_{(p)}$ and $T_{(p)}^*$
satisfy the assumptions of Corollary~\ref{fin_rep_seq},
both (i) and (ii),
for $p$ and $p'$, respectively.
\qed

\medskip
\noindent
Institute of Mathematics, Technical University,
Wroc{\l}aw, Poland

\noindent and

\noindent
Department of Mathematics, University of Alberta,
Edmonton, Alberta, Canada T6G 2G1.


\end{document}